\documentclass[12pt,fleqn]{article}
\usepackage{mathrsfs}
\usepackage{amsfonts}
\usepackage{amssymb}
\usepackage{latexsym,amsmath}
\usepackage{graphicx}
\date{}
\begin{document}
\title{On the normalized Laplacian spectra of some subdivision joins of two graphs\footnote{Corresponding author.  E-mail: guixiantian@163.com (G.-X.
Tian), cuishuyu@zjnu.cn (S.-Y. Cui)}}
\author{Gui-Xian Tian$^a$, Shu-Yu Cui$^b$,\\
{\small{\it $^a$College of Mathematics, Physics and Information Engineering,}}\\
{\small{\it Zhejiang Normal University, Jinhua, Zhejiang, 321004,
P.R. China}} \\
{\small{\it $^b$Xingzhi College, Zhejiang Normal University, Jinhua,
Zhejiang, 321004, P.R. China}} }\maketitle

\begin{abstract} For two simple graphs $G_1$ and $G_2$, we denote
the subdivision-vertex join and subdivision-edge join of $G_1$ and
$G_2$ by $G_1\dot{\vee}G_2$ and $G_1\veebar G_2$, respectively. This
paper determines the normalized Laplacian spectra of
$G_1\dot{\vee}G_2$ and $G_1\veebar G_2$ in terms of these of $G_1$
and $G_2$ whenever $G_1$ and $G_2$ are regular. As applications, we
construct some non-regular normalized Laplacian cospectral graphs.
Besides we also compute the number of spanning trees and the
degree-Kirchhoff index of $G_1\dot{\vee}G_2$ and $G_1\veebar G_2$
for regular graphs $G_1$ and $G_2$.



\end{abstract}

\section*{1. Introduction}
\indent\indent Throughout this paper, all graphs considered are
finite simple graphs. Let $G=(V,E)$ be a graph with vertex set
$V=\{v_{1},v_{2},\ldots,v_{n}\}$ and edge set $E(G)$. The adjacency
matrix $A(G)$ of $G$ is an $n\times n$ matrix whose $(i,j)$-entry is
$1$ if $v_{i}$ and $v_{j}$ are adjacent in $G$ and $0$ otherwise.
The degree of $v_{i}$ in $G$ is denoted by $d_{i}=d_{G}(v_{i})$. Let
$D(G)$ be the degree diagonal matrix of $G$ with diagonal entries
$d_{1},d_{2},\ldots,d_{n}$. The \emph{normalized Laplacian matrix}
$\mathcal {L}(G)$ of $G$ is defined as
$I_n-D(G)^{-1/2}A(G)D(G)^{-1/2}$, where $I_n$ denotes the identity
matrix of order $n$. Denote the characteristic polynomial
$\det(xI_n-\mathcal {L}(G))$ of $\mathcal {L}(G)$ by $\phi(G;x)$.
The roots of $\phi(G;x)$ are called the \emph{normalized Laplacian
eigenvalues} of $G$, denoted by $0=\lambda_1 (G)\leq \lambda_2
(G)\leq\cdots\leq \lambda_n (G)\leq2$. The set of all eigenvalues is
called the \emph{normalized Laplacian spectrum} of $G$.

The normalized Laplacian matrix $\mathcal {L}(G)$, which is
consistent with the transition probability matrix
$P(G)=D(G)^{-1}A(G)$ in the random walk on $G$ and spectral
geometry\cite{Chung1997} has attracted people's attention. For
example, Banerjee and Jost\cite{Banerjee2008} studied how the
normalized Laplacian spectrum is affected by operations such as
motif doubling, graph splitting and joining. Butler and
Grout\cite{Butler2011} constructed many pairs of non-regular
normalized Laplacian cospectral graphs. Cavers et
al.\cite{Cavers2010} obtained some bounds on the normalized
Laplacian energy and general Randi\'{c} index $R_{-1}$ of graphs.
Chen et al.\cite{Chen2004} gave an interlacing inequality on the
normalized Laplacian eigenvalues of $G$. Chen and
Zhang\cite{Chen2007} obtained the formulae for the resistance
distance and degree-Kirchoff index in terms of the normalized
Laplacian eigenvalues and eigenvectors of $G$.
Kirkland\cite{Kirkland2006} investigated the Limit points for
normalized Laplacian eigenvalues and so on. For more review about
the normalized Laplacian spectrum of graphs, readers may refer to
\cite{Chung1997}. Recently, Chen and Liao\cite{Chen2017} determined
the normalized Laplacian spectra of the (edge)corona for two graphs.
Almost at the same time, the normalized Laplacian spectra of some
subdivision-coronas for two regular graphs were determined by Das
and Panigrahi\cite{Das2017}. The \emph{subdivision
graph}\cite{Cvetkovic2010} $S(G)$ of a graph $G$ is the graph
obtained by inserting a new vertex into every edge of $G$. Denote
the set of such new vertices by $I(G)$. The following graph
operations based on subdivision graphs come from
\cite{Indulal2012}.\\
\\
\textbf{Definition 1.1} The \emph{subdivision-vertex join}
$G_1\dot{\vee}G_2$ of two graphs $G_1$ and $G_2$ is the graph
obtained from $S(G_1)$ and $G_2$ by joining each vertex of $V(G_1)$
with every vertex of $V(G_2)$.\\
\\
\textbf{Definition 1.2} The \emph{subdivision-edge join} $G_1\veebar
G_2$  of two graphs $G_1$ and $G_2$ is the graph obtained from
$S(G_1)$ and $G_2$ by joining each vertex of $I(G_1)$ with every
vertex of $V(G_2)$.

In \cite{Indulal2012}, the spectra of $G_1\dot{\vee}G_2$ and
$G_1\veebar G_2$ were computed in terms of these of regular graphs
$G_1$ and $G_2$. The author also constructed infinite family of new
integral graphs. Bu et al.\cite{Bu2014} obtained the formulae for
resistance distance of $G_1\dot{\vee}G_2$ and $G_1\veebar G_2$.
Recently, Liu and Zhang\cite{Liu2017} determined the spectra,
(signless) Laplacian spectra of $G_1\dot{\vee}G_2$ and $G_1\veebar
G_2$ for a regular graph $G_1$ and arbitrary graph $G_2$. As
applications, they constructed infinitely many pairs of cospectral
graphs and obtained the number of spanning trees and the Kirchhoff
index of $G_1\dot{\vee}G_2$ and $G_1\veebar G_2$.

Motivated by the works above, we focus on determining the normalized
Laplacian spectra of $G_1\dot{\vee}G_2$ and $G_1\veebar G_2$ in
terms of those of regular graphs $G_1$ and $G_2$(see Theorems 2.3
and 2.4, respectively). As applications, we construct some
non-regular normalized Laplacian cospectral graphs(see Theorem 3.1
and Example 3.2). Besides we also compute the number of spanning
trees and the degree-Kirchhoff index of $G_1\dot{\vee}G_2$ and
$G_1\veebar G_2$ for regular graphs $G_1$ and $G_2$(see Theorems 3.3
and 3.4, respectively).

\section*{2. Main results}
\indent\indent In this section, we shall determine the normalized
Laplacian spectra of $G_1\dot{\vee}G_2$ and $G_1\veebar G_2$ with
the help of the $coronal$ of matrices. Recall that the $M$-coronal
$\Gamma_{M}(x)$\cite{Cui2012,McLeman2011} of a matrix $M$ of order
$n$ is defined to be the sum of the entries of the matrix
$(xI_{n}-M)^{-1}$, that is, ${\Gamma _M}(x) = 1_n^T{(x{I_n} - M)^{ -
1}}{1_n}$, where $1_{n}$ denotes the column vector of size $n$ with
all entries equal to one. For the sake of convenience, we first give
a slight generalization of the $M$-coronal of a matrix $M$ of order
$n$ below:
\[
\Gamma _{M} (x,\alpha) = 1_n^T (xI_n  - M \circ (\alpha J_n  + (1 -
\alpha )I_n ))^{ - 1} 1_n,
\]
where $J_n$ denotes the matrix of order $n$ with all entries equal
to one and ``$\circ$" the Hadamard product\cite{Horn} of two
matrices. If $\alpha=1$, then $\Gamma _{M} (x,1)$ is the usual
$M$-coronal $\Gamma_{M}(x)$ introduced by Cui and Tian in
\cite{Cui2012}. Let $G$ be an $r$-regular graph of order $n$. Then
\begin{equation}\label{1}
\Gamma _{\mathcal {L}(G)} (x,\alpha ) = {n \mathord{\left/
 {\vphantom {n {(x + \alpha  - 1)}}} \right.
 \kern-\nulldelimiterspace} {(x + \alpha  - 1)}}.
\end{equation}

Let $G=(V,E)$ be a simple graph with vertex set
$V=\{v_{1},v_{2},\ldots,v_{n}\}$ and edge set
$E(G)=\{e_1,e_2,\ldots,e_m\}$. The \emph{incident matrix} of $G$ is
an $n\times m$ matrix whose $(i,j)$-entry is  $1$ if $v_{i}$ and
$e_{j}$ are incident in $G$ and $0$ otherwise, denoted by $R(G)$. If
the graph $G$ is an $r$ regular, then $R(G)R(G)^{T}=A(G)+rI_n$. The
\emph{line graph} of a graph $G$ is the graph $l(G)$, whose vertices
are the edges of $G$ and two vertices of $l(G)$ are adjacent if and
only if they are incident on a common vertex in $G$. It is well
known\cite{Cvetkovic2010} that $R(G)^{T}R(G)=A(l(G))+2I_m$.

The following Lemmas 2.1 and 2.2 used in the rest of article come
from \cite{Das2017} and \cite{Liu2017}, respectively.\\
\\
\textbf{Lemma 2.1} {\it Let $G$ be an $r$-regular graph on $n$
vertices and $m$ edges. Then the eigenvalues of $A(l(G))$ are the
eigenvalues of $2(r-1)I_n-r\mathcal {L}(G)$ and $-2$ repeated $m-n$
times.}\\
\\
\textbf{Lemma 2.2} {\it Let $A$ be a matrix of order $n$. Then
\[
\det (xI_n  - A - \alpha J_n ) = (1 - \alpha \Gamma _A (x,1))\det
(xI_n  - A).
\]}

We first consider the normalized Laplacian spectrum of
$G_1\dot{\vee}G_2$. Let $G_1$ be an $r_{1}$-regular graph on $n_1$
vertices and $m_{1}$ edges, and $G_2$ be an $r_{2}$-regular graph on
$n_2$ vertices. Then
\[
A(G_1 \dot  \vee G_2 ) = \left( {\begin{array}{*{20}c}
   {0_{n_1 } } & R & {J_{n_1  \times n_2 } }  \\
   {R^T } & {0_{m_1 } } & {0}  \\
   {J_{n_2  \times n_1 } } & {0 } & {A(G_2 )}  \\
\end{array}} \right),
\]
and
\[
D(G_1 \dot  \vee G_2 ) = \left( {\begin{array}{*{20}c}
   {(r_1  + n_2 )I_{n_1 } } & 0 & 0  \\
   0 & {2I_{m_1 } } & 0  \\
   0 & 0 & {(r_2  + n_1 )I_{n_2 } }  \\
\end{array}} \right).
\]
Thus we have
\begin{equation}\label{2}
\mathcal {L}(G_1 \dot  \vee G_2 ) = \left( {\begin{array}{*{20}c}
   {I_{n_1 } } & { - \frac{R}{{\sqrt {2(r_1  + n_2 )} }}} & { - \frac{{J_{n_1  \times n_2 } }}{{\sqrt {(r_1  + n_2 )(r_2  + n_1 )} }}}  \\
   { - \frac{{R^T }}{{\sqrt {2(r_1  + n_2 )} }}} & {I_{m_1 } } & 0  \\
   { - \frac{{J_{n_2  \times n_1 } }}{{\sqrt {(r_1  + n_2 )(r_2  + n_1 )} }}} & 0 & {\mathcal {L}(G_2 ) \circ B} \\
\end{array}} \right),
\end{equation}
where $B = \alpha J_{n_2 }  + (1 - \alpha )I_{n_2 }$ with
$\alpha=r_2/(r_2+n_1)$.\\
\\
\textbf{Theorem 2.3} {\it Let $G_1$ be an $r_{1}$-regular graph on
$n_1$ vertices and $m_{1}$ edges, and $G_2$ be an $r_{2}$-regular
graph on $n_2$ vertices. Also let $0=\mu_1,\mu_2,\cdots,\mu_{n_1}$
and $0=\nu_1,\nu_2,\cdots,\nu_{n_2}$ be the normalized Laplacian
spectrum of $G_1$ and $G_2$, respectively. Then the normalized
Laplacian spectra of $G_1\dot{\vee}G_2$ consists of:

(i) $0$;

(ii) $1$, repeated $m_{1}-n_{1}$ times;

(iii) $ \frac{{n_1 +r_2 \nu _i }}{{r_2  + n_1 }} $, for
$i=2,\ldots,n_{2}$;

(iv) two roots of the equation, for $i=2,\ldots,n_{1}$,
\[
(x - 1)^2  - \frac{{r_1 (2 - \mu _i )}}{{2(r_1  + n_2 )}} = 0;
\]

(v) two roots of the equation
\[
x^2  - (2 + \frac{{n_1 }}{{r_2  + n_1 }})x + \frac{{2n_1 }}{{r_2  +
n_1 }} + \frac{{n_2 r_2 }}{{(r_1  + n_2 )(r_2  + n_1 )}} = 0.
\]
}
\\
\textbf{Proof.}  According to (\ref{2}), we obtain the
characteristic polynomial of $\mathcal {L}(G_1 \dot  \vee G_2 )$ is
\begin{equation}\label{3}
\begin{array}{l}
 \phi (G_1 \dot  \vee G_2 ;x) = \det (xI_{n_1  + m_1  + n_2 }  - \mathcal {L}(G_1 \dot  \vee G_2 ))
 \\\\
\;\;\; \;\;\;\;\;\;\;\;\;\;\;\;\;\;\;\;\;\; = \det \left(
{\begin{array}{*{20}c}
   {(x - 1)I_{n_1 } } & {\frac{R}{{\sqrt {2(r_1  + n_2 )} }}} & {\frac{{J_{n_1  \times n_2 } }}{{\sqrt {(r_1  + n_2 )(r_2  + n_1 )} }}}  \\
   {\frac{{R^T }}{{\sqrt {2(r_1  + n_2 )} }}} & {(x - 1)I_{m_1 } } & 0  \\
   {\frac{{J_{n_2  \times n_1 } }}{{\sqrt {(r_1  + n_2 )(r_2  + n_1 )} }}} & 0 & {xI_{n_2 }  - \mathcal {L}(G_2 ) \circ B}  \\
\end{array}} \right) \\\\
\;\;\; \;\;\;\;\;\;\;\;\;\;\;\;\;\;\;\;\;\; = \det (xI_{n_2 }  - \mathcal {L}(G_2 ) \circ B)\cdot\det (S), \\
 \end{array}
\end{equation}
where
\[
\begin{array}{l}
 S = \left( {\begin{array}{*{20}c}
   {(x - 1)I_{n_1 }  - \frac{{J_{n_1  \times n_2 } }}{{\sqrt {(r_1  + n_2 )(r_2  + n_1 )} }}(xI_{n_2 }  - \mathcal {L}(G_2 ) \circ B)^{ - 1} \frac{{J_{n_2  \times n_1 } }}{{\sqrt {(r_1  + n_2 )(r_2  + n_1 )} }}} & {\frac{R}{{\sqrt {2(r_1  + n_2 )} }}}  \\
   {\frac{{R^T }}{{\sqrt {2(r_1  + n_2 )} }}} & {(x - 1)I_{m_1 } }  \\
\end{array}} \right) \\\\
\;\;\;\; = \left( {\begin{array}{*{20}c}
   {(x - 1)I_{n_1 }  - \frac{{\Gamma _{\mathcal {L}(G_2 )} \left( {x,{{r_2 } \mathord{\left/
 {\vphantom {{r_2 } {(r_2  + n_1 )}}} \right.
 \kern-\nulldelimiterspace} {(r_2  + n_1 )}}} \right)}}{{(r_1  + n_2 )(r_2  + n_1 )}}J_{n_1 } } & {\frac{R}{{\sqrt {2(r_1  + n_2 )} }}}  \\
   {\frac{{R^T }}{{\sqrt {2(r_1  + n_2 )} }}} & {(x - 1)I_{m_1 } }  \\
\end{array}} \right) \\
 \end{array}
\]
is the Schur complement\cite{Zhang} of $ {xI_{n_2 }  - \mathcal
{L}(G_2 ) \circ B}$. Notice that $RR^{T}=2r_1I_{n_1}- {r_1\mathcal
{L}(G_1 )}$ as $G_1$ is $r_1$-regular. Then, from Lemma 2.2 and
(\ref{1}), one has
\begin{equation}\label{4}
\begin{array}{l}
 \det (S) = \det \left( {\begin{array}{*{20}c}
   {(x - 1)I_{n_1 }  - \frac{{\Gamma _{\mathcal {L}(G_2 )} \left( {x,{{r_2 } \mathord{\left/
 {\vphantom {{r_2 } {(r_2  + n_1 )}}} \right.
 \kern-\nulldelimiterspace} {(r_2  + n_1 )}}} \right)}}{{(r_1  + n_2 )(r_2  + n_1 )}}J_{n_1 } } & {\frac{R}{{\sqrt {2(r_1  + n_2 )} }}}  \\
   {\frac{{R^T }}{{\sqrt {2(r_1  + n_2 )} }}} & {(x - 1)I_{m_1 } }  \\
\end{array}} \right) \\\\
 \;\;\;\;\;\;\;\;\;\;\;= (x - 1)^{m_1 } \det \left( {(x - 1)I_{n_1 }  - \frac{{\Gamma _{\mathcal {L}(G_2 )} \left( {x,{{r_2 } \mathord{\left/
 {\vphantom {{r_2 } {(r_2  + n_1 )}}} \right.
 \kern-\nulldelimiterspace} {(r_2  + n_1 )}}} \right)}}{{(r_1  + n_2 )(r_2  + n_1 )}}J_{n_1 }  - \frac{{RR^T }}{{2(r_1  + n_2 )(x - 1)}}} \right)
 \\\\
  \;\;\;\;\;\;\;\;\;\;\; = (x - 1)^{m_1 } \det \left( {(x - 1 - \frac{r_1}{{(r_1  + n_2 )(x - 1)}})I_{n_1 }  + \frac{{r_1\mathcal {L}(G_1 )}}{{2(r_1  + n_2 )(x - 1)}} - \frac{{\Gamma _{\mathcal {L}(G_2 )} \left( {x,{{r_2 } \mathord{\left/
 {\vphantom {{r_2 } {(r_2  + n_1 )}}} \right.
 \kern-\nulldelimiterspace} {(r_2  + n_1 )}}} \right)}}{{(r_1  + n_2 )(r_2  + n_1 )}}J_{n_1 } } \right)
 \\\\
   \;\;\;\;\;\;\;\;\;\;\;= (x - 1)^{m_1 } \left( {1 - \frac{{\Gamma _{\mathcal {L}(G_2 )} \left( {x,{{r_2 } \mathord{\left/
 {\vphantom {{r_2 } {(r_2  + n_1 )}}} \right.
 \kern-\nulldelimiterspace} {(r_2  + n_1 )}}} \right)}}{{(r_1  + n_2 )(r_2  + n_1 )}}\Gamma _{\mathcal {L}(G_1 )} (x - 1 - \frac{r_1}{{(r_1  + n_2 )(x - 1)}},1)}
 \right)\\\\
  \;\;\;\;\;\;\;\;\;\;\;\;\;\;\;\;\;\;\times\det \left( {(x - 1 - \frac{r_1}{{(r_1  + n_2 )(x - 1)}})I_{n_1 }  + \frac{{r_1\mathcal {L}(G_1 )}}{{2(r_1  + n_2 )(x - 1)}}} \right) .\\
 \end{array}
\end{equation}
Again, from (\ref{1}), one has
\begin{equation}\label{5}
\Gamma _{\mathcal {L}(G_1 )} (x - 1 - \frac{{r_1 }}{{(r_1  + n_2 )(x
- 1)}},1) = \frac{{n_1 }}{{x - 1 - \frac{{r_1 }}{{(r_1  + n_2 )(x -
1)}}}}
\end{equation}
and
\begin{equation}\label{6}
\Gamma _{\mathcal {L}(G_2 )} \left( {x,\frac{{r_2 }}{{r_2  + n_1 }}}
\right) = \frac{{n_2 }}{{x - 1 + \frac{{r_2 }}{{r_2  + n_1 }}}}.
\end{equation}
Now, substituting (\ref{5}) and (\ref{6}) into (\ref{4}), we obtain
\begin{equation}\label{7}
\begin{array}{l}
 \det (S) = (x - 1)^{m_1  - n_1 } \left( {1 - \frac{{n_1 n_2 (x - 1)}}{{(r_1  + n_2 )(r_2  + n_1 )\left( {x - \frac{{n_1 }}{{r_2  + n_1 }}} \right)\left( {(x - 1)^2  - \frac{{r_1 }}{{r_1  + n_2 }}} \right)}}} \right)
 \\\\
 \;\;\;\;\;\;\;\;\;\;\;\;\;\;\;\;\;\times\det \left( {\left( {(x - 1)^2  - \frac{{r_1 }}{{r_1  + n_2 }}} \right)I_{n_1 }  + \frac{{r_1 }}{{2(r_1  + n_2 )}}\mathcal {L}(G_1 )} \right) .\\
 \end{array}
\end{equation}
On the other hand, it is easy to see that
\begin{equation}\label{8}
\det (xI_{n_2 }  - \mathcal {L}(G_2 ) \circ B) = \det \left( {\left(
{x - \frac{{n_1 }}{{r_2  + n_1 }}} \right)I_{n_2 }  - \frac{{r_2
}}{{r_2 + n_1 }}\mathcal {L}(G_2 )} \right).
\end{equation}
Thus, substituting (\ref{7}) and (\ref{8}) into (\ref{3}), we obtain
the characteristic polynomial of $\mathcal {L}(G_1 \dot  \vee G_2 )$
is
\[
\begin{array}{l}
 \phi (G_1 \dot  \vee G_2 ;x) = (x - 1)^{m_1  - n_1 } \left( {\left( {x - \frac{{n_1 }}{{r_2  + n_1 }}} \right)\left( {(x - 1)^2  - \frac{{r_1 }}{{r_1  + n_2 }}} \right) - \frac{{n_1 n_2 (x - 1)}}{{(r_1  + n_2 )(r_2  + n_1 )}}} \right)
 \\\\
 \;\;\;\;\;\;\;\;\;\;\;\;\;\;\;\;\;\;\;\;\;\;\;\;\;\;\;\;\; \times\prod\limits_{i = 2}^{n_1 } {\left( {(x - 1)^2  - \frac{{r_1 }}{{r_1  + n_2 }} + \frac{{r_1 \mu _i }}{{2(r_1  + n_2 )}}} \right)} \prod\limits_{i = 2}^{n_2 } {\left( {x - \frac{{n_1 }}{{r_2  + n_1 }} - \frac{{r_2 \nu _i }}{{r_2  + n_1 }}} \right)}
 \\\\
 \;\;\;\;\;\;\;\;\;\;\;\;\;\;\;\;\;\;\;\;= x(x - 1)^{m_1  - n_1 } \left( {x^2  - (2 + \frac{{n_1 }}{{r_2  + n_1 }})x + \frac{{2n_1 }}{{r_2  + n_1 }} + \frac{{n_2 r_2 }}{{(r_1  + n_2 )(r_2  + n_1 )}}} \right)
  \\\\
\;\;\;\;\;\;\;\;\;\;\;\;\;\;\;\;\;\;\;\;\;\;\;\;\;\;\;\;\;\times
\prod\limits_{i = 2}^{n_1 } {\left( {(x - 1)^2  - \frac{{r_1(2- \mu
_i) }}{{2(r_1  + n_2 )}}} \right)} \prod\limits_{i = 2}^{n_2 }
{\left( {x - \frac{{n_1+r_2 \nu _i }}{{r_2 + n_1 }}} \right)},
 \\
 \end{array}
\]
which implies the required results. $\Box$\\

Next we consider the normalized Laplacian spectrum of $G_1\veebar
G_2$. Let $G_1$ be an $r_{1}$-regular graph on $n_1$ vertices and
$m_{1}$ edges, and $G_2$ be an $r_{2}$-regular graph on $n_2$
vertices. Then
\[
A(G_1 \veebar G_2 ) = \left( {\begin{array}{*{20}c}
   {0_{n_1 } } & R & 0  \\
   {R^T } & {0_{m_1 } } & {J_{m_1  \times n_2 } }  \\
   0 & {J_{n_2  \times m_1 } } & {A(G_2 )}  \\
\end{array}} \right)
\]
and
\[
D(G_1\veebar G_2 ) = \left( {\begin{array}{*{20}c}
   {r_1 I_{n_1 } } & 0 & 0  \\
   0 & {(2 + n_2 )I_{m_1 } } & 0  \\
   0 & 0 & {(r_2  + m_1 )I_{n_2 } }  \\
\end{array}} \right).
\]
Thus we have
\begin{equation}\label{9}
L(G_1 \veebar G_2 ) = \left( {\begin{array}{*{20}c}
   {I_{n_1 } } & { - \frac{R}{{\sqrt {r_1 (2 + n_2 )} }}} & 0  \\
   { - \frac{{R^T }}{{\sqrt {r_1 (2 + n_2 )} }}} & {I_{m_1 } } & { - \frac{{J_{m_1  \times n_2 } }}{{\sqrt {(2 + n_2 )(r_2  + m_1 )} }}}  \\
   0 & { - \frac{{J_{n_2  \times m_1 } }}{{\sqrt {(2 + n_2 )(r_2  + m_1 )} }}} & {\mathcal {L}(G_2 ) \circ B}  \\
\end{array}} \right)
\end{equation}
where $B = \alpha J_{n_2 }  + (1 - \alpha )I_{n_2 }$ with
$\alpha=r_2/(r_2+m_1)$.\\
\\
\textbf{Theorem 2.4} {\it Let $G_1$ be an $r_{1}$-regular graph on
$n_1$ vertices and $m_{1}$ edges, and $G_2$ be an $r_{2}$-regular
graph on $n_2$ vertices. Also let $0=\mu_1,\mu_2,\cdots,\mu_{n_1}$
and $0=\nu_1,\nu_2,\cdots,\nu_{n_2}$ be the normalized Laplacian
spectra of $G_1$ and $G_2$, respectively. Then the normalized
Laplacian spectrum of $G_1\veebar G_2$ consists of:

(i) $0$;

(ii) $1$, repeated $m_{1}-n_{1}$ times;

(iii) $ {\frac{{m_1  + r_2 \nu _i }}{{r_2  + m_1 }}}$, for
$i=2,\ldots,n_{2}$;

(iv) two roots of the equation, for $i=2,\ldots,n_{1}$,
\[
(x - 1)^2  - \frac{{2 - \mu _i }}{{2 + n_2 }} = 0;
\]

(v) two roots of the equation
\[
x^2  - (2 + \frac{{m_1 }}{{r_2  + m_1 }})x + \frac{{2m_1 }}{{r_2  +
m_1 }} + \frac{{n_2 r_2 }}{{(2  + n_2 )(r_2  + m_1 )}} = 0.
\]}
\\
\textbf{Proof.}  By (\ref{9}), the characteristic polynomial of
$\mathcal {L}(G_1 \veebar G_2 )$ is
\begin{equation}\label{10}
\begin{array}{l}
 \phi (G_1 \veebar  G_2 ;x) = \det (xI_{n_1  + m_1  + n_2 }  -\mathcal {L}(G_1 \veebar  G_2 ))
 \\\\
 \;\;\;\;\;\;\;\;\;\;\;\;\;\;\;\;\;\;\;\;\;\; = \det \left( {\begin{array}{*{20}c}
   {(x - 1)I_{n_1 } } & {\frac{R}{{\sqrt {r_1 (2 + n_2 )} }}} & 0  \\
   {\frac{{R^T }}{{\sqrt {r_1 (2 + n_2 )} }}} & {(x - 1)I_{m_1 } } & {\frac{{J_{m_1  \times n_2 } }}{{\sqrt {(2 + n_2 )(r_2  + m_1 )} }}}  \\
   0 & {\frac{{J_{n_2  \times m_1 } }}{{\sqrt {(2 + n_2 )(r_2  + m_1 )} }}} & {xI_{n_2 }  - \mathcal {L}(G_2 ) \circ B}  \\
\end{array}} \right) \\\\
 \;\;\;\;\;\;\;\;\;\;\;\;\;\;\;\;\;\;\;\;\;\; = \det (xI_{n_2 }  - \mathcal {L}(G_2 ) \circ B)\det (S), \\
 \end{array}
\end{equation}
where
\[
\begin{array}{l}
 S = \left( {\begin{array}{*{20}c}
   {(x - 1)I_{n_1 } } & {\frac{R}{{\sqrt {r_1 (2 + n_2 )} }}}  \\
   {\frac{{R^T }}{{\sqrt {r_1 (2 + n_2 )} }}} & {(x - 1)I_{m_1 }  - \frac{{J_{m_1  \times n_2 } }}{{\sqrt {(2 + n_2 )(r_2  + m_1 )} }}(xI_{n_2 }  - \mathcal {L}(G_2 ) \circ B)^{ - 1} \frac{{J_{n_2  \times m_1 } }}{{\sqrt {(2 + n_2 )(r_2  + m_1 )} }}}  \\
\end{array}} \right) \\\\
 \;\;\;\; = \left( {\begin{array}{*{20}c}
   {(x - 1)I_{n_1 } } & {\frac{R}{{\sqrt {r_1 (2 + n_2 )} }}}  \\
   {\frac{{R^T }}{{\sqrt {r_1 (2 + n_2 )} }}} & {(x - 1)I_{m_1 }  - \frac{{\Gamma _{\mathcal {L}(G_2 )} \left( {x,{{r_2 } \mathord{\left/
 {\vphantom {{r_2 } {(r_2  + m_1 )}}} \right.
 \kern-\nulldelimiterspace} {(r_2  + m_1 )}}} \right)}}{{(2 + n_2 )(r_2  + m_1 )}}J_{m_1 } }  \\
\end{array}} \right) \\
 \end{array}
\]
is the Schur complement\cite{Zhang} of $ {xI_{n_2 }  - \mathcal
{L}(G_2 ) \circ B}$. Then, from Lemma 2.2 and (\ref{1}), one has
\begin{equation}\label{11}
\begin{array}{l}
 \det (S) = \det \left( {\begin{array}{*{20}c}
   {(x - 1)I_{n_1 } } & {\frac{R}{{\sqrt {r_1 (2 + n_2 )} }}}  \\
   {\frac{{R^T }}{{\sqrt {r_1 (2 + n_2 )} }}} & {(x - 1)I_{m_1 }  - \frac{{\Gamma _{\mathcal {L}(G_2 )} \left( {x,{{r_2 } \mathord{\left/
 {\vphantom {{r_2 } {(r_2  + m_1 )}}} \right.
 \kern-\nulldelimiterspace} {(r_2  + m_1 )}}} \right)}}{{(2 + n_2 )(r_2  + m_1 )}}J_{m_1 } }  \\
\end{array}} \right) \\\\
  \;\;\;\;\;\;\;\;\;\;\;  = (x - 1)^{n_1 } \det \left( {(x - 1)I_{m_1 }  - \frac{{R^T R}}{{r_1 (2 + n_2 )(x - 1)}} - \frac{{\Gamma _{\mathcal {L}(G_2 )} \left( {x,{{r_2 } \mathord{\left/
 {\vphantom {{r_2 } {(r_2  + m_1 )}}} \right.
 \kern-\nulldelimiterspace} {(r_2  + m_1 )}}} \right)}}{{(2 + n_2 )(r_2  + m_1 )}}J_{m_1 } } \right)
 \\\\
  \;\;\;\;\;\;\;\;\;\;\;  = (x - 1)^{n_1 } \left( {1 - \frac{{\Gamma _{\mathcal {L}(G_2 )} \left( {x,{{r_2 } \mathord{\left/
 {\vphantom {{r_2 } {(r_2  + m_1 )}}} \right.
 \kern-\nulldelimiterspace} {(r_2  + m_1 )}}} \right)}}{{(2 + n_2 )(r_2  + m_1 )}}\Gamma _{\frac{{R^T R}}{{r_1 (2 + n_2 )(x - 1)}}} (x - 1,1)} \right)
 \\\\
  \;\;\;\;\;\;\;\;\;\;\; \;\;\;\;\;\;\times\det \left( {(x - 1)I_{m_1 }  - \frac{{R^T R}}{{r_1 (2 + n_2 )(x - 1)}}} \right) .\\
 \end{array}
\end{equation}
It is follows from (\ref{1}) that
\begin{equation}\label{12}
\Gamma _{\mathcal {L}(G_2 )} \left( {x,\frac{{r_2 }}{{r_2  + m_1 }}}
\right) = \frac{{n_2 }}{{x - 1 + \frac{{r_2 }}{{r_2  + m_1 }}}}
\end{equation}
and
\begin{equation}\label{13}
\Gamma _{\frac{{R^T R}}{{r_1 (2 + n_2 )(x - 1)}}} (x - 1,1) =
\frac{{m_1 }}{{x - 1 - \frac{{2r_1 }}{{r_1 (2 + n_2 )(x - 1)}}}}.
\end{equation}
Notice that $R^{T}R=A(l(G_1))+2I_{m_1}$. Now, substituting
(\ref{12}) and (\ref{13}) into (\ref{11}), along with Lemma 2.1, we
obtain
\begin{equation}\label{14}
\begin{array}{l}
 \det (S) = (x - 1)^{n_1 } \left( {1 - \frac{{m_1 n_2 (x - 1)}}{{(2 + n_2 )(r_2  + m_1 )\left( {x - 1 + \frac{{r_2 }}{{r_2  + m_1 }}} \right)\left( {(x - 1)^2  - \frac{2}{{2 + n_2 }}} \right)}}} \right)  \\
  \;\;\;\;\;\;\;\;\;\;\; \;\;\;\;\;\;\;\times\prod\limits_{i = 1}^{m_1 } {\left( {x - 1 - \frac{{2 + \tau _i (l(G_1 ))}}{{r_1 (2 + n_2 )(x - 1)}}} \right)}
  \\
   \;\;\;\;\;\;\;\;\;\;\;= (x - 1)^{m_1  - n_1 } \left( {1 - \frac{{m_1 n_2 (x - 1)}}{{(2 + n_2 )(r_2  + m_1 )\left( {x -  \frac{{m_1 }}{{r_2  + m_1 }}} \right)\left( {(x - 1)^2  - \frac{2}{{2 + n_2 }}} \right)}}} \right)
   \\
 \;\;\;\;\;\;\;\;\;\;\;\;\;\;\;\;\;\; \times \prod\limits_{i = 2}^{n_1 } {\left( {(x - 1)^2  - \frac{{2 - \mu _i }}{{2 + n_2 }}} \right)} \left( {(x - 1)^2  - \frac{2}{{2 + n_2 }}} \right) .\\
 \end{array}
\end{equation}
where $\tau_i(l(G_1))$ denotes an eigenvalue of the line graph
$l(G_1)$ of $G_1$.

On the other hand, by a simple computation, one gets
\begin{equation}\label{15}
\det (xI_{n_2 }  - \mathcal {L}(G_2 ) \circ B) = \det \left( {\left(
{x - \frac{{m_1 }}{{r_2  + m_1 }}} \right)I_{n_2 }  - \frac{{r_2
}}{{r_2 + m_1 }}\mathcal {L}(G_2 )} \right).
\end{equation}
Now, substituting (\ref{14}) and (\ref{15}) into (\ref{10}), we
obtain the characteristic polynomial of $\mathcal {L}(G_1\dot \vee
G_2)$ is
\[
\begin{array}{l}
 \phi (G_1 \underline  \vee  G_2 ;x) = (x - 1)^{m_1  - n_1 } \left( {\left( {x - \frac{{m_1 }}{{r_2  + m_1 }}} \right)\left( {(x - 1)^2  - \frac{2}{{2 + n_2 }}} \right) - \frac{{m_1 n_2 (x - 1)}}{{(2 + n_2 )(r_2  + m_1 )}}} \right)
 \\\\
 \;\;\;\;\;\;\;\;\;\;\;\;\;\;\;\;\;\;\;\;\;\;\;\;\;\;\;\;\;\;\;\;\times\prod\limits_{i = 2}^{n_1 } {\left( {(x - 1)^2  - \frac{{2 - \mu _i }}{{2 + n_2 }}} \right)} \prod\limits_{i = 2}^{n_2 } {\left( {x - \frac{{m_1 }}{{r_2  + m_1 }} - \frac{{r_2 \nu _i }}{{r_2  + m_1 }}} \right)}
 \\\\
 \;\;\;\;\;\;\;\;\;\;\;\;\;\;\;\;\;\;\;\; = x(x - 1)^{m_1  - n_1 } \left( {x^2  - (2 + \frac{{m_1 }}{{r_2  + m_1 }})x + \frac{{2m_1 }}{{r_2  + m_1 }} + \frac{{n_2 r_2 }}{{(2  + n_2 )(r_2  + m_1 )}}} \right)
 \\\\
 \;\;\;\;\;\;\;\;\;\;\;\;\;\;\;\;\;\;\;\;\;\;\;\;\;\;\;\;\;\;\;\;\times\prod\limits_{i = 2}^{n_1 } {\left( {(x - 1)^2  - \frac{{2 - \mu _i }}{{2 + n_2 }}} \right)} \prod\limits_{i = 2}^{n_2 } {\left( {x -\frac{{m_1+r_2 \nu _i }}{{r_2  + m_1 }}}
 \right)},
 \\
 \end{array}
\]
which follows the desired results. $\Box$

\section*{3. Applications}

\indent\indent In this section, we consider to construct many pairs
of non-regular normalized Laplacian cospectral graphs. Besides, we
also compute the number of spanning trees and the degree-Kirchhoff
index of $G_1\dot{\vee}G_2$ and $G_1\veebar G_2$ for regular graphs
$G_1$ and $G_2$.

Two graphs $G$ and $H$ are called the \emph{normalized Laplacian
cospectral} if the normalized Laplacian spectra of $G$ and $H$ are
the same. It is well known that many families of pairs of cospectral
graphs may be constructed by using some graph operations(for
example, see\cite{Butler2011,Cui2012,Das2017,Liu2017,McLeman2011}).
Recall that two regular graphs are normalized Laplacian cospectral
if and only if they are cospectral. Hence what's more interesting is
to construct non-regular normalized Laplacian cospectral graphs. The
following Theorem 3.1 may construct many pairs of non-regular
normalized
Laplacian cospectral graphs. From Theorems 2.3 and 2.4, we arrive immediately at:\\
\\
\textbf{Theorem 3.1} {\it Let $G_1$ and $G_2$ (not necessarily
distinct) be two normalized Laplacian cospectral regular graphs.
Also let $H_1$ and $H_2$ (not necessarily distinct) be two
normalized Laplacian cospectral regular graphs, then
$G_1\dot{\vee}H_1$ and $G_2\dot{\vee}H_2$ (respectively, $G_1\veebar
H_1$ and $G_2\veebar
H_2$) are normalized Laplacian cospectral.}\\
\\
\textbf{Example 3.2} Let $Q_4$ and $K_2$ be the $4$-cube (tesseract)
and complete graph of order $2$. Also let $\widehat{Q_4}$(drawn in
[2, p.15]) be the graph obtained from $Q_4$ by the Godsil-McKay
switching. Since $Q_4$ and $\widehat{Q_4}$ are nonisomorphic
$4$-regular cospectral graphs\cite{Brouwer2011}. Thus they are also
normalized Laplacian cospectral. Notice that both $Q_4\dot{\vee}
K_2$ and $\widehat{Q_4}\dot{\vee} K_2$ are non-regular graphs of
order $50$ with the degree sequence $(6^{(16)},2^{(32)},17^{(2)})$,
where $a^{(b)}$ indicates that $a$ is repeated $b$ times. Applying
Theorem 3.1, we obtain $Q_4\dot{\vee} K_2$ and
$\widehat{Q_4}\dot{\vee} K_2$ are non-regular normalized Laplacian
cospectral graphs. Similarly, $Q_4\veebar K_2$ and
$\widehat{Q_4}\veebar K_2$ are also non-regular normalized Laplacian
cospectral graphs with the degree sequence $(4^{(48)},33^{(2)})$.\\

Let $G$ be a connected graph on $n$ vertices and $m$ edges, and the
normalized Laplacian spectrum $0=\lambda_1 (G)<\lambda_2
(G)\leq\cdots\leq \lambda_n (G)$. It is proved\cite{Chen2007} that
the degree-Kirchhoff index of $G$ is
\begin{equation}\label{16}
Kf^*(G) = 2m\sum\limits_{i = 2}^{n } {\frac{1}{{\lambda_i(G) }}}.
\end{equation}
Recently, the degree Kirchhoff index of some graph operations have
been investigated, such as line graphs, subdivision graph, total
graphs, coronae, edge coronae graphs and so on(for example,
see\cite{Chen2017,Huang2015}). Next we shall give the formulae for
the degree-Kirchhoff index of the vertex-join $G_1\dot{\vee}G_2$ and
edge-join $G_1\veebar G_2$ of regular graphs $G_1$ and $G_2$.\\
\\
\textbf{Theorem 3.3} {\it Let $G_1$ be an $r_{1}$-regular graph on
$n_1$ vertices and $m_{1}$ edges, and $G_2$ be an $r_{2}$-regular
graph on $n_2$ vertices. Also let $0=\mu_1,\mu_2,\cdots,\mu_{n_1}$
and $0=\nu_1,\nu_2,\cdots,\nu_{n_2}$ be the normalized Laplacian
spectra of $G_1$ and $G_2$, respectively. Then
\[
\begin{array}{l}
 Kf^*(G_1\dot{\vee}G_2) = (n_1 r_1  + n_2 r_2  + 2n_1 n_2  + 2m_1 )\\
  \; \;\;\;\;\;\;\;\; \;\;\;\;\;\;\;\;\; \;\;\;\;\;\;\;\times\left( {m_1  - n_1  + \frac{{(r_1  + n_2 )(2r_2  + 3n_1 )}}{{2n_1 r_1  + 2n_1 n_2  + n_2 r_2 }}  + \sum\limits_{i = 2}^{n_1 } {\frac{{4(r_1  + n_2 )}}{{2n_2  + r_1 \mu _i }}}+ \sum\limits_{i = 2}^{n_2 } {\frac{{r_2  + n_1 }}{{n_1  + r_2 \nu _i }}}  } \right), \\
 \end{array}
\]
and
\[
\begin{array}{l}
 Kf^*(G_1\veebar G_2) = (n_1 r_1  + n_2 r_2  + 2m_1 n_2  + 2m_1 )\\
\; \;\;\;\;\;\;\;\; \;\;\;\;\;\;\;\;\; \;\;\;\;\;\;\;\;\; \times\left( {m_1  - n_1  + \frac{{(2 + n_2 )(2r_2  + 3m_1 )}}{{4m_1  + 2m_1 n_2  + n_2 r_2 }} + \sum\limits_{i = 2}^{n_1 } {\frac{{4 + 2n_2 }}{{n_2  + \mu _i }}}  + \sum\limits_{i = 2}^{n_2 } {\frac{{r_2  + m_1 }}{{m_1  + r_2 \nu _i }}} } \right). \\
 \end{array}
\]}
\textbf{Proof.} By the Vieta theorem, two roots
$x_1^{(i)},x_2^{(i)}$ of the equation
\[
(x - 1)^2  - \frac{{r_1 (2 - \mu _i )}}{{2(r_1  + n_2 )}} =
0,\;\;\;\;(i=2,\ldots,n_{1})
\]
satisfy
\[
x_1^{(i)} x_2^{(i)}  = 1 - \frac{{r_1 (2 - \mu _i )}}{{2(r_1  + n_2
)}},\;\;x_1^{(i)}  + x_2^{(i)}  = 2.
\]
Similarly, two roots $x_1,x_2$ of the equation
\[
x^2  - (2 + \frac{{n_1 }}{{r_2  + n_1 }})x + \frac{{2n_1 }}{{r_2  +
n_1 }} + \frac{{n_2 r_2 }}{{(r_1  + n_2 )(r_2  + n_1 )}} = 0
\]
satisfy
\[
x_1 x_2  = \frac{{2n_1 }}{{r_2  + n_1 }} + \frac{{n_2 r_2 }}{{(r_1 +
n_2 )(r_2  + n_1 )}},\;\;x_1  + x_2  = 2 + \frac{{n_1 }}{{r_2  + n_1
}}.
\]
Thus, from Theorem 2.3 and (\ref{16}), we obtain
\[
\begin{array}{l}
 Kf^*(G_1\dot{\vee}G_2) = (n_1 (r_1  + n_2 ) + 2m_1  + n_2 (r_2  + n_1 ))\sum\limits_{i = 2}^{n_1  + m_1  + n_2 } {\frac{1}{{\lambda _i (G_1\dot{\vee}G_2)}}}
 \\\\
  \; \;\;\;\;\;\;\;\; \;\;\;\;\;\;\;\;\;\;\;\;= (n_1 (r_1  + n_2 ) + 2m_1  + n_2 (r_2  + n_1 ))\\\\
  \; \;\;\;\;\;\;\;\; \;\;\;\;\;\;\;\;\; \;\;\;\;\;\;\;\times\left( {m_1  - n_1   + \frac{{x_1  + x_2 }}{{x_1 x_2 }}+ \sum\limits_{i = 2}^{n_1 } {\frac{{x_1^{(i)}  + x_2^{(i)} }}{{x_1^{(i)} x_2^{(i)} }} + } \sum\limits_{i = 2}^{n_2 } {\frac{{r_2  + n_1 }}{{n_1  + r_2 \nu _i }}} } \right)
  \\\\
  \; \;\;\;\;\;\;\;\; \;\;\;\;\;\;\;\;\;\;\;\;= (n_1 r_1  + n_2 r_2  + 2n_1 n_2  + 2m_1 )\\\\
  \; \;\;\;\;\;\;\;\; \;\;\;\;\;\;\;\;\; \;\;\;\;\;\;\;\times\left( {m_1  - n_1  + \frac{{(r_1  + n_2 )(2r_2  + 3n_1 )}}{{2n_1 r_1  + 2n_1 n_2  + n_2 r_2 }}  + \sum\limits_{i = 2}^{n_1 } {\frac{{4(r_1  + n_2 )}}{{2n_2  + r_1 \mu _i }}}+ \sum\limits_{i = 2}^{n_2 } {\frac{{r_2  + n_1 }}{{n_1  + r_2 \nu _i }}}  } \right). \\
 \end{array}
\]
Using the same technique, the formulae of $Kf^*(G_1\veebar G_2)$
follows from
Theorem 2.4 and (\ref{16}). $\Box$\\

Let $G$ be a connected graph on $n$ vertices and $m$ edges, and the
normalized Laplacian spectra $0=\lambda_1 (G)<\lambda_2
(G)\leq\cdots\leq \lambda_n (G)\leq2$. It is well
known\cite{Chung1997} that the number of spanning trees of $G$ is
\begin{equation}\label{17}
t(G) = \frac{1}{{2m}}\prod\limits_{i = 1}^n {d_i } \prod\limits_{j =
2}^n {\lambda _j }.
\end{equation}
The number of spanning trees of graph operations has been studied
extensively, such as line graphs, subdivision graph, total graphs,
coronae, edge coronae graphs and so on(for example,
see\cite{Chen2017,Huang2015,Liu2017}). Now we give the number of
spanning trees of the vertex-join $G_1\dot{\vee}G_2$ and edge-join
$G_1\veebar G_2$ for regular graphs $G_1$ and $G_2$. Remark that The
number of spanning trees of $G_1\dot{\vee}G_2$ and $G_1\veebar G_2$
has also been computed by using their Laplacian
spectra\cite{Liu2017}.
\\
\\
\textbf{Theorem 3.4} {\it Let $G_1$ be an $r_{1}$-regular graph on
$n_1$ vertices and $m_{1}$ edges, and $G_2$ be an $r_{2}$-regular
graph on $n_2$ vertices. Also let $0=\mu_1,\mu_2,\cdots,\mu_{n_1}$
and $0=\nu_1,\nu_2,\cdots,\nu_{n_2}$ be the normalized Laplacian
spectra of $G_1$ and $G_2$, respectively. Then
\[
t(G_1\dot{\vee}G_2) = 2^{m_1  - n_1  + 1} \prod\limits_{i = 2}^{n_1 } {(2n_2  + r_1 \mu _i )} \prod\limits_{j = 2}^{n_2 } {(n_1  + r_2 \nu _i )},\\
\]
and
\[
t(G_1\veebar G_2) = r_1 ^{n_1 } (2 + n_2 )^{m_1  - n_1 } \prod\limits_{i = 2}^{n_1 } {(n_2  + \mu _i )} \prod\limits_{j = 2}^{n_2 } {(m_1  + r_2 \nu _i )} . \\
\]
}\\
\\
\textbf{Proof.} From Theorem 2.3, the Vieta theorem and (\ref{17}),
we get easily
\[
\begin{array}{l}
 t(G_1\dot{\vee}G_2) = \frac{{2^{m_1 } (r_1  + n_2 )^{n_1 } (r_2  + n_1 )^{n_2 } \left( {\frac{{2n_1 }}{{r_2  + n_1 }} + \frac{{n_2 r_2 }}{{(r_1  + n_2 )(r_2  + n_1 )}}}
 \right)\prod\nolimits_{i = 2}^{n_1 } {\left( {1 - \frac{{r_1 (2 - \mu _i )}}{{2(r_1  + n_2 )}}} \right)} \prod\nolimits_{j = 2}^{n_2 } {\left( {\frac{{n_1  + r_2 \nu _i }}{{r_2  + n_1 }}} \right)} }}{{n_1 (r_1  + n_2 ) + 2m_1  + n_2 (r_2  + n_1 )}} \\\\
  \;\;\;\;\;\;\;\; \;\;\;\;\;\;\;\;= 2^{m_1  - n_1  + 1} \prod\limits_{i = 2}^{n_1 } {(2n_2  + r_1 \mu _i )} \prod\limits_{j = 2}^{n_2 } {(n_1  + r_2 \nu _i )} . \\
 \end{array}
\]
Similarly, it follows from Theorem 2.4 and and (\ref{17}) that
\[
\begin{array}{l}
 t(G_1\veebar G_2) = \frac{{r_1 ^{n_1 } (2 + n_2 )^{m_1 } (r_2  + m_1 )^{n_2 } \left( {\frac{{2m_1 }}{{r_2  + m_1 }} + \frac{{n_2 r_2 }}{{(2 + n_2 )(r_2  + m_1 )}}} \right)
 \prod\nolimits_{i = 2}^{n_1 } {\left( {1 - \frac{{2 - \mu _i }}{{2 + n_2 }}} \right)} \prod\nolimits_{j = 2}^{n_2 } {\left( {\frac{{m_1  + r_2 \nu _i }}{{r_2  + m_1 }}} \right)} }}{{r_1 n_1  + (2 + n_2 )m_1  + n_2 (r_2  + m_1 )}}
 \\\\
 \; \;\;\;\;\;\;\;\; \;\;\;\;\;\;\;\;= r_1 ^{n_1 } (2 + n_2 )^{m_1  - n_1 } \prod\limits_{i = 2}^{n_1 } {(n_2  + \mu _i )} \prod\limits_{j = 2}^{n_2 } {(m_1  + r_2 \nu _i )} . \\
 \end{array}
\]
Hence, the result follows. $\Box$

\section*{4. Conclusion}

\indent\indent We remark that all nonzero off-diagonal elements of
the normalized Laplacian matrix are not integer. It is not easy to
determine the normalized Laplacian spectrum of some graph operations
in terms of those of factor graphs. In this paper, we first
introduce a slight generalization of the $M$-coronal of a matrix $M$
of order $n$
\[
\Gamma _{M} (x,\alpha) = 1_n^T (xI_n  - M \circ (\alpha J_n  + (1 -
\alpha )I_n ))^{ - 1} 1_n.
\]
This conception is used to determine the normalized Laplacian
spectra of $G_1\dot{\vee}G_2$ and $G_1\veebar G_2$ in terms of those
of factor graphs $G_1$ and $G_2$. As applications, some non-regular
normalized Laplacian cospectral graphs are constructed. Finally, we
also compute the number of spanning trees and the degree-Kirchhoff
index of $G_1\dot{\vee}G_2$ and $G_1\veebar G_2$ for two regular
graphs $G_1$ and $G_2$.

And more notably, using the generalized $M$-coronal technique here,
we may discuss to describe the normalized Laplacian spectra of some
graph operations, such as join,(edge)corona, neighbourhood corona,
variants of (edge)corona, variants of neighbourhood (edge)corona,
variants of join based on Q-graph(R-graph, total graph, duplication
graph) and so on. Remark that the case of (edge)corona has been
studied by describing all corresponding eigenvectors with their
eigenvalues in \cite{Chen2017}.
\\
\\
\textbf{Acknowledgements} This work was in part supported by NNSFC
(Nos. 11371328, 11671053) and by the Natural Science Foundation of
Zhejiang Province, China (No. LY15A010011).


{\small \begin{thebibliography}{16}

\bibitem{Banerjee2008} A. Banerjee, J. Jost, On the spectrum of the normalized graph
Laplacian, Linear Algebra Appl. 428 (2008) 3015-3022.

\bibitem{Brouwer2011} A.E. Brouwer, W.H. Haemers, Spectra of Graphs, Springer, New York,
2011.

\bibitem{Bu2014} C.-J. Bu, B. Yan, X.-Q. Zhou, J. Zhou, Resistance distance in
subdivision-vertex join and subdivision-edge join of graphs, Linear
Algebra Appl. 458 (2014) 454-462.

\bibitem{Butler2011} S. Butler, J. Grout, A construction of cospectral graphs for the
normalized Laplacian, Electron. J. Combin. 18 (1) (2011)
$\sharp$P231.

\bibitem{Cavers2010} M.S. Cavers, S. Fallat, S. Kirkland, On the normalized Laplacian
energy and general Randi\'{c} index $R_{-1}$ of graphs, Linear
Algebra Appl. 433 (2010) 172-190.

\bibitem{Chen2004} G.T. Chen, G. Davis, F. Hall, Z.S. Li, K. Patel, M. Stewart, An
interlacing result on normalized Laplacians, SIAM J. Discrete Math.
18(2) (2004) 353-361.

\bibitem{Chen2017} H.Y. Chen, L.W. Liao, The normalized Laplacian spectra of
the corona and edge corona of two graphs, Linear Multilinear
Algebra, 65 (2017) 582-592.

\bibitem{Chen2007} H.Y. Chen, F.J. Zhang, Resistance distance and the normalized
Laplacian spectrum, Discrete Appl. Math. 155 (2007) 654-661.

\bibitem{Chung1997} F.R.K. Chung, Spectral Graph Theory, CBMS Regional Conference Series in Mathematics, Amer. Math. Soc., Providence, 1997.

\bibitem{Cui2012} S.-Y. Cui, G.-X. Tian, The spectrum and the signless
Laplacian spectrum of coronae, Linear Algebra Appl., 437 (2012)
1692-1703.

\bibitem{Cvetkovic2010} D. M. Cvetkovi\'{c}, P. Rowinson, H. Simi\'{c}, An
introduction to the Theory of Graph Spectra, Cambridge University
Press, Cambridge, 2010.

\bibitem{Das2017} A. Das, P. Panigrahi, Normalized Laplacian spectrum
of some subdivision-coronas of two regular graphs, Linear
Multilinear Algebra 65 (2017) 962-972.

\bibitem{Horn} R. A. Horn, C. R. Johnson, Topics in matrix analysis, Cambridge University Press, 1991.

\bibitem{Huang2015} J. Huang, S. Li, On the normalised Laplacian spectrum,
degree-Kirchhoff index and spanning trees of graphs, Bull. Aust.
Math. Soc. 91 (2015) 353-367.

\bibitem{Indulal2012} G. Indulal, Spectrum of two new joins of graphs and infinite families of integral graphs, Kragujevac J. Math. 36
(2012) 133-139.

\bibitem{Kirkland2006} S. Kirkland, Limit points for normalized Laplacian eigenvalues,
Electron. J. Linear Algebra 15 (2006) 337-344.

\bibitem{Liu2017} X.-G. Liu, Z.H. Zhang, Spectra of subdivision-vertex join and subdivision-edge join of two graphs, Bull. Malays. Math. Sci. Soc. (2017).
doi:10.1007/s40840-017-0466-z.

\bibitem{McLeman2011} C. McLeman, E. McNicholas, Spectra of coronae, Linear Algebra
Appl. 435 (2011) 998-1007.

\bibitem{Zhang} F.-Z. Zhang, The Schur complement and its applications, Springer, 2005.





\end{thebibliography}}
\end{document}